\documentclass[twocolumn]{autart}
\usepackage[longnamesfirst,compress]{natbib}
\usepackage{latexsym, amssymb, amsbsy, amsopn, amstext}
\usepackage{amsmath}
\usepackage{graphicx, subfigure}
\usepackage{pgf,tikz}
\usepackage{picins}
\usetikzlibrary{arrows,shapes,automata,positioning}
\usetikzlibrary{fit}
\usepackage{pst-node}
\newtheorem{assumption}{Assumption}
\newtheorem{lemma}{Lemma}
\newtheorem{theorem}{Theorem}
\newtheorem{remark}{Remark}

\begin{document}

\begin{frontmatter}

\title{Distributed Output Regulation for a Class of Nonlinear Multi-Agent Systems with Unknown-Input Leaders\thanksref{footnoteinfo}}

\thanks[footnoteinfo]{This work was supported in part by the Fundamental Research Funds for the Central Universities under Grant 24820152015RC36 and in part by National Natural Science Foundation of China under Grants 61503033, 61333001 and 61503359. This paper was not presented at any IFAC meeting. Tel.: +86-10-62651449. Fax.: +86-10-62587343.}

\author[bupt]{Yutao Tang}\ead{yttang@bupt.edu.cn},    
\author[cas]{Yiguang Hong}\ead{yghong@iss.ac.cn},               
\author[ustc]{Xinghu Wang}\ead{xinghuw@ustc.edu.cn}  
\address[bupt]{School of Automation, Beijing University of Posts and Telecommunications, Beijing 100876, China}
\address[cas]{Key Laboratory of Systems and Control, Institute of Systems Science, Chinese Academy of Sciences, Beijing 100190, China}
\address[ustc]{Department of Automation, University of Science and Technology of China, Hefei 230027, China}

\begin{keyword}
Multi-agent systems, distributed output regulation, nonlinear dynamics, internal model
\end{keyword}

\begin{abstract}
In this paper, a distributed output regulation problem is formulated for a class of uncertain nonlinear multi-agent systems subject to local disturbances. The formulation is given to study a leader-following problem when the leader contains unknown inputs and its dynamics is different from those of the followers.  Based on the conventional output regulation assumptions and graph theory, distributed feedback controllers are constructed to make the agents globally or semi-globally follow the uncertain leader even when the bound of the leader's inputs is unknown to the followers.
\end{abstract}

\end{frontmatter}

\section{Introduction}

The past decade has witnessed a rapid development in the field of multi-agent systems and fruitful results have been achieved for the leader-following coordination problem. In recent years, distributed output regulation of multi-agent systems has been proposed to provide a general framework for leader-following consensus in linear or nonlinear cases \citep{hong2013rnc, dong2014asme}. Internal model approach, developed to solve the conventional output regulation \citep{isidori2003, huang2004}, was effectively used for the distributed design, especially for nonlinear agents. For example, a cooperative output regulation problem was considered for a class of nonlinear uncertain multi-agent systems with unity relative degree in \citet{su2013scl}, while agents in the output-feedback form were considered with an autonomous leader and no-loop graph \citep{ding2013tac}.

It may be restrictive or unpractical if we always consider an autonomous leader
without unknown inputs, especially in the case when the leader is an
uncooperative target or contains unmodeled uncertainties.
Therefore, it is necessary to study multi-agent control when the
leader contains (unknown) inputs.  In fact, the
generalized output regulation (GOR) problem to track an exosystem (or a
leader) with external inputs was discussed in many
publications including \citet{saberi2001ijc} and \citet{ramos2004tac}. On the other hand, a distributed problem was investigated when the agent dynamics are double
integrators to track a leader with an unknown but bounded acceleration
in \citet{cao2012tac}, and then a similar design was given in
\citet{li2013tac} by assuming that the leader and followers share the
same dynamics. To our best of knowledge, there are no general
results on nonlinear multi-agent control when the unknown-input leader and the followers have different dynamics with external disturbances.

The objective of our paper is to study distributed output regulation
for leader-following multi-agent systems with an unknown-input leader,
whose dynamics is nonlinear and may differ from those of the
followers. The contribution of the work is at least twofold:
\begin{itemize}
\item We extend the distributed output regulation to the case when the leader contains unknown inputs and has a dynamics different from those of the non-identical followers with (unbounded) local disturbances, and provide distributed controls to solve this problem in different cases. The results are consistent with the existing output regulation results when the leader does not have unknown inputs and the disturbances are bounded \citep[e.g.][]{dong2014asme,su2013scl}.
\item We extend the conventional GOR to its distributed version for multi-agent systems with an unknown-input leader. Moreover, both global and semi-global results are obtained for nonlinear agents with unity relative degree, while only local results were obtained for a conventional (single-agent) case in \citet{ramos2004tac}.
\end{itemize}

Notations: Let $\mathbb{R}^n$ be the $n$-dimensional Euclidian
space and $\mathbb{R}^{n}_M=\{s\in \mathbb{R}^{n}\mid -M\leq
s_i\leq M, i=1,{\dots},n\}$ for a constant $M>0$. For a vector
$x$, $||x||$ (or $||x||_\infty$) denotes its Euclidian norm (or
infinite norm). $\text{diag}\{b_1,{\dots},b_n\}$ denotes an $n\times n$
diagonal matrix with diagonal elements $b_i\; (i=1,{\dots},n)$;
$\text{col}(a_1,{\dots},a_n) = [a_1^T,{\dots},a_n^T]^T$ for column
vectors $a_i\; (i=1,{\dots},n)$.  A continuous function $\alpha\colon[0,\, a)\to [0,\, \infty)$ belongs to
class $\mathcal{K}$ if it is strictly increasing and $\alpha(0)=0$;
It belongs to class $\mathcal{K}_\infty$ if it belongs to class $\mathcal{K}$ with $a=\infty$ and $\lim_{s\to \infty}\alpha(s)\to\infty$.

\section{Problem Formulation}

Consider a group of $n+1$ agents with one leader (regarded as node 0) as follows:
\begin{equation}\label{exosystem}
\dot{v}=p(v)+q(v)w(t),\quad y_0=r(v,\mu)
\end{equation}
where $v\in \mathbb{R}^{n_v}$ is the leader's state, and $w(t)\in \mathbb{R}^{n_w}, y_0\in \mathbb{R}$ are its input and output, respectively.  Here $\mu\in \mathbb{R}^{n_\mu}$ is an uncertain parameter vector, and $w(t)$ is continuous satisfying $||w(t)||_{\infty}\leq l$ with a constant $l>0$. The other $n$ (non-identical) agents are followers described by
\begin{align}\label{follower}
\begin{cases}
  \dot{z}_i=f_i(z_i,y_i,\mu)\\
  \dot{y}_i=g_i(z_i,y_i,\mu)+ b_iu_i+d_i, \qquad i=1,{\dots},n
  \end{cases}
\end{align}
where $z_i\in \mathbb{R}^{n_{z_i}},\; y_i \in \mathbb{R},\; d_i\in \mathbb{R},\;
b_i>0$. Without loss of generality, we take $b_i=1$ and assume all functions $f_i$, $g_i$, $p$, $q$, $r$ are smooth with $f_i(0,0,\mu)=0$, $g_i(0,0,\mu)=0$, $p(0)=0$, $r(0,\mu)=0$. $d_i$ is the local disturbance of agent $i$ governed by
\begin{equation}\label{disturbance}
\dot{\omega}_i=S_i\omega_i,\quad d_i=D_i(\mu)\omega_i.
\end{equation}
As usual, we assume $S_i\in \mathbb{R}^{n_{\omega_i}\times n_{\omega_i}}$ has no eigenvalues with negative real parts \citep{huang2004}.

Clearly, the first-order nonlinear agent in \citet{liu2013consensus} is a special case of \eqref{follower}, and system (\ref{follower}) was also considered to solve an output consensus problem for the exosystem without any inputs in \citet{ding2013tac}.

The interaction topology among these agents can be described by a graph $\mathcal{G}=(\mathcal{V}, \mathcal{E})$, where $\mathcal{V}=\{0,1,{\dots},n\}$ is the set of nodes and $\mathcal{E}$ is the set of arcs. $(i,j)$ denotes an arc leaving from node $i$ and entering node $j$ \citep{godsil2001}. A walk in
graph $\mathcal {G}$ is an alternating sequence $i_{1}e_{1}i_{2}e_{2}{\cdots}e_{k-1}i_{k}$ of nodes $i_{l}$ and arcs $e_{m}=(i_{m},i_{m+1})\in\mathcal {E}$ for $l=1,2,{\dots},k$. If there exists a walk from node $i$ to node $j$ then node $i$ is said to be reachable from $j$.  Define the neighbor set of agent $i$ as $\mathcal{N}_i=\{j: (j,i)\in \mathcal {E} \}$ for $i=1,\dots,n$.  A weighted adjacency matrix of $\mathcal {G}$ is denoted by $A=[a_{ij}]\in \mathbb{R}^{(n+1)\times (n+1)}$, where $a_{ii}=0$ and $a_{ij}\geq 0$($a_{ij}>0$ if $(j,i)\in \mathcal{E}$).  A graph is said to be undirected if $a_{ij}=a_{ji}$ ($i,j=0,1,{\dots},n$). The Laplacian $L=[l_{ij}]\in \mathbb{R}^{(n+1)\times (n+1)}$ of graph  $\mathcal{G}$ is defined as $l_{ii}=\sum_{j\neq i}a_{ij}$ and $l_{ij}=-a_{ij} (j\neq i)$.  Denote $\bar{\mathcal{G}}$ as the induced subgraph of $\mathcal{G}$ associated with the $n$ followers. The following assumption has been widely used in coordination of multi-agent systems \citep{hong2006tracking, su2013scl}.

\begin{assumption} \label{ass1:graph} The leader (node 0) is reachable from any other node of
$\mathcal{G}$ and the induced subgraph $\bar{\mathcal G}$ is undirected.
\end{assumption}

Given a communication graph $\mathcal{G}$, denote $H\in \mathbb{R}^{n\times n}$ as the submatrix of the Laplacian $L$ by deleting its first row and first column. By Lemma 3 in \citet{hong2006tracking}, $H$ is positive definite under Assumption \ref{ass1:graph}.  The distributed control law can be constructed as follows:
\begin{equation}\label{ctr-nominal}
\begin{split}
  u_i&=k_i(\xi_i, y_i-y_j, j\in \mathcal{N}_i),\\
  \dot{\xi}_i&=h_i(\xi_i, y_i-y_j, j\in \mathcal{N}_i)
\end{split}
\end{equation}
where $\xi_i\in \mathbb{R}^{n_{\xi_i}}$ with a nonnegative integer $n_{\xi_i}$ and functions
$k_i(\cdot), h_i(\cdot)$ to be designed later.

To handle this nonlinear multi-agent system with an unknown-input leader, we formulate the problem as the {\it distributed generalized output regulation problem} or simply {\it distributed regulation problem}.  It is said to be (globally) solved for systems (\ref{exosystem})-(\ref{disturbance}) with a given graph $\mathcal{G}$, if we can find a distributed control law (\ref{ctr-nominal}), such that, for any $(z_i(0), y_i(0))\in \mathbb{R}^{n_{z_i}+1}$, $\mu\in \mathbb{R}^{n_{\mu}}$, $\xi_i(0)\in \mathbb{R}^{n_{\xi_i}}$, $v(0)\in \mathbb{R}^{n_v}, \omega_i(0)\in \mathbb{R}^{n_{\omega_i}}$, the trajectory of the closed-loop system, composed of (\ref{exosystem}), (\ref{follower}), and (\ref{ctr-nominal}), is well-defined for all $t>0$, and moreover,
\begin{equation} \label{errori}
\lim_{t\to \infty}e_i(t)=0,\quad e_i=y_i-y_0, \quad i=1,{\dots},n.
\end{equation}
Our problem is said to be semi-globally solved if, for any given $M>0$, we can find a control law
(\ref{ctr-nominal}) with a constant $\bar M\geq 0$, such that, for any initial condition
$(z_i(0), y_i(0))\in \mathbb{R}_M^{n_{z_i}+1}$, $\mu\in\mathbb{R}^{n_{\mu}}_M$, $\xi_i(0)\in \mathbb{R}_{\bar M}^{n_{{\xi}_i}}$, $v(0)\in \mathbb{R}_M^{n_v}$ and $\omega_i(0)\in \mathbb{R}_M^{n_{\omega_i}}$, the trajectory of the closed-loop system is well-defined for all $t>0$ and (\ref{errori}) holds.

\begin{remark}
When $n=1$, our problem becomes GOR studied in \citet{saberi2001ijc} and \citet{ramos2004tac}. Here we seek non-local output feedback control for nonlinear systems of the form \eqref{follower}, while only local results were obtained in \citet{ramos2004tac} requiring the exosystem's state.
\end{remark}
\begin{remark}
Because of the unknown expression or type of $w(t)$, adaptive IM discussed in \citet{su2013scl} fails to solve our problem even with a time-varying one  like in \citet{yang2010cdc}. In fact, our problem can be viewed as a distributed version of GOR to handle those exosystems (or leaders) with unknown inputs, which certainly extends the existing multi-agent output regulation formulation when the leaders have no unknown inputs \citep{hong2013rnc,su2013scl}.
\end{remark}

The following assumption was used for GOR of nonlinear systems \citep[see][]{ramos2004tac}.

\begin{assumption}\label{ass3:exo}
There exist two class $\mathcal{K}$ functions $\alpha_0(\cdot)$ and $\gamma_0(\cdot)$ such that
$$||v(t)||\leq \alpha_0(||v(0)||)+\gamma_0(||w(t)||_\infty).$$
\end{assumption}
Clearly, $v(t)$ is bounded by $\bar b_1=\alpha_0(||v(0)||)+\gamma_0(l)$ and $\bar b_2=\max_{v\in
\mathbb{R}^{n_v}_{\bar b_1}}\{\frac{\partial {r}}{\partial
v}q(v)\}$ is well-defined under this assumption. Moreover, system \eqref{exosystem} is Lyapunov stable at $v=0$ when $w=0$ with the neutrally stable exosystem for nonlinear output regulation as one of its special cases.

Similar to the output regulation problem, the solvability of regulator equations plays a key role in the study of nonlinear GOR. Therefore, we give the following assumption for the solution of regulator equations.

\begin{assumption}\label{ass4:re}
For $i=1,{\dots},n$, there exists a smooth function ${\bf z }_i(v, \mu)$ with ${\bf z}_i(0, \mu)=0$ such that,
\begin{equation}\label{eq:re}
\frac{\partial {\bf z}_i(v, \mu)}{\partial v}p(v)=f_i({\bf z}_i(v, \mu), r(v, \mu)).
\end{equation}
\end{assumption}

Under Assumption \ref{ass4:re},  letting ${\bf u}_i(v, \mu)=\scalebox{1}{$\frac{\partial  {r}(v, \mu)}{\partial v}$}p(v)-g_i({\bf z}_i(v, \mu), r(v, \mu))$ and performing a coordinate transformation:
 $\bar z_i=z_i-{\bf z}_i(v, \mu)$,\; $e_i=y_i-r(v, \mu)$ gives
\begin{align}\label{err}
\begin{cases}
  \dot{\bar z}_i=\bar f_i(\bar z_i, e_i, v, w, \mu)\\
  \dot{e}_i=\bar g_i(\bar z_i,e_i, v, w, \mu)+u_i+d_i
\end{cases}
\end{align}
where
\begin{align*}
\!\bar f_i(\bar z_i, e_i, v, w, \mu)&= \hat f_i(\bar z_i,e_i,v, \mu)-
\scalebox{1}{$\frac{\partial  {\bf z}_i(v, \mu)}{\partial v}$}q(v)w\\
\!\bar g_i(\bar z_i, e_i, v, w, \mu)&= \hat g_i(\bar z_i,e_i,v, \mu)-{\bf u}_i(v, \mu)-\scalebox{1}{$\frac{\partial  r(v, \mu)}{\partial  v}$}q(v)w\\
  \hat f_i(\bar z_i, e_i,v, \mu)&=f_i(z_i,y_i, \mu)-f_i({\bf z}_i(v, \mu), r(v, \mu), \mu)\\
  \hat g_i(\bar z_i, e_i,v, \mu)&=g_i(z_i,y_i, \mu)-g_i({\bf z}_i(v, \mu),{r}(v, \mu), \mu).
\end{align*}
Furthermore, $\bar f_i(0,0, v, 0, \mu)=0$, $\bar g_i(0,0,v,0, \mu)=0$, $\hat g_i(0,0,v, \mu)=0$. For simplicity, the error system can be rewritten as
\begin{equation}\label{blockform}
\begin{cases}
\dot{\bar z}=\bar f(\bar z,e,v,w, \mu)\\
\dot{e}=\bar g(\bar z,e,v,w, \mu)+ u+d
\end{cases}
\end{equation}
where $\bar z=\text{col}(\bar z_1,{\dots},\bar z_n)$, $e=\text{col}(e_1,{\dots},e_n)$, $u=\text{col}(u_1,{\dots}, u_n)$, $d=\text{col}(d_1, \dots, d_n)$ and $\bar f$, $\bar g$ are suitably defined by (\ref{err}).

By (\ref{blockform}), output regulation problem of the multi-agent system (\ref{exosystem}), (\ref{follower}) and (\ref{disturbance}) is transformed into a problem of finding a distributed control law in the form of (\ref{ctr-nominal}) such that $\lim_{t\to\infty} {e(t)} =0$ and the trajectory of the closed-loop system is well-defined for $t>0$. For this purpose, we introduce the following assumption for the zero dynamics (i.e., the $\bar z$-subsystem) of system (\ref{blockform}), though we need not ``stabilize" this subsystem (to make $\bar z$ vanish).

\begin{assumption}\label{ass6:unknown}
For any compact subset $\Sigma\subset \mathbb{R}^{\hat n}$ $(\hat n=n_v+n_w+n_\mu)$, and for $i=1,{\dots},n$, there exists a smooth Lyapunov function $V_{\bar z_i}(\cdot)$ satisfying
$\alpha_{1i}(||\bar z_i||)\leq V_{\bar z_i}(\bar z_i)\leq
\alpha_{2i}(||\bar z_i||)$ for some smooth functions $\alpha_{1i}(\cdot)$, $\alpha_{2i}(\cdot)\in \mathcal{K}_\infty$, such that, for any $(v(t),w(t),\mu)\in \Sigma$
  \begin{equation}\label{eq:unknown}
  \dot{V}_{\bar z_i}(\bar z_i)\mid_{(\ref{err})}\leq -\alpha_i(||\bar z_i||)+\delta_{1i}\gamma_{1i}(e_i)+\delta_{2i}\gamma_{2i}(w)
  \end{equation}
where $\gamma_{1i}(\cdot)$, $\gamma_{2i}(\cdot)$ are known smooth positive definite functions, $\alpha_i(\cdot)$ is a known class $\mathcal{K}_\infty$ function satisfying $\limsup_{s \to 0+}(\alpha_i^{-1}(s^2)/s)\leq \infty$, and $\delta_{1i},\delta_{2i}$ are some unknown positive constants.
\end{assumption}

\begin{remark}\label{rem0}
Although the $\bar z_i$-subsystem is related with $v$, the
condition (\ref{eq:unknown}) is not restrictive since $v(t)$ is bounded by Assumption
\ref{ass3:exo} and $\bar f_i(0,0, v, 0,\mu)=0$.  Similar assumptions were commonly used in the study of nonlinear output regulation \cite*[e.g.][]{dong2014asme, su2013scl, xu2010tcs}.
\end{remark}

\section{Main Results}

In this section, we give a constructive design to solve our problem by investigating
system (\ref{blockform}).

For the following non-smooth analysis, consider an equation $\dot{x}=f(x,t)$ with a discontinuous righthand side, where $f\colon \mathbb{R}^m\times \mathbb{R}\to \mathbb{R}^m$ is measurable and essentially locally bounded. By Proposition 3 in \citet{cortes2008ieeecs}, it has a Filippov solution on $[t_0,\, t_1]$.  Let $V\colon \mathbb{R}^m\to \mathbb{R}$ be a locally Lipschitz continuous function. $ \dot{V}\triangleq \bigcap_{\xi\in \partial V}\xi^T\mathcal{F}[f](z,t)$ represents the set-valued Lie derivative of $V$, where $\partial V(z)$ denotes the Clarke's generalized gradient of $V$ \citep{cortes2008ieeecs}.

Then we start with a general discussion and give a simpler design for a special case.
\subsection{General Discussion}

Here each agent only knows that the unknown leader's input is bounded, but does not know the exact value of the input bound.

To track the leader and meanwhile reject the (unbounded) local disturbance $d_i$, we split the total control effort into two parts as $u_i=u_i^d+u_i^r$, where the term $u_i^d$ is to deal with $d_i$ and $u_i^r$ to make $y_i$ follow $y_0$. It is well-known that internal model methods were effective to reject modeled disturbances. Here we construct $u_i^d$ following the same technical line.

Let $P_i(s)=s^{n_{p_i}}+\hat{p}_{i,1}s^{n_{p_i}-1}+\cdots+\hat{p}_{i,n_{p_i}-1}s+\hat{p}_{i,n_{p_i}}$ be the minimal polynomial of the matrix $S_i$ and denote $\tau_i=\text{col}(\tau_{i,1},\ldots,\tau_{i,n_{p_i}})$ with $\tau_{i,j}=\frac{{\rm d}^{j-1}d_i(t)}{{\rm d}t^{j-1}}$.  Take
\begin{align*}
  \Phi_i=\left[\begin{array}{c|c}
    0&I_{n_{p_i}-1}\\\hline
    -\hat{p}_{i,n_{p_i}}&\mathbf{\hat p}_i
  \end{array}\right],\;
  \Psi_i=\begin{bmatrix}
    1&0_{1\times (n_{p_i}-1)}
  \end{bmatrix}
\end{align*}
where $\mathbf{\hat p}_i=[\,-\hat{p}_{i,n_{p_i}-1}\,\cdots\,-\hat{p}_{i,1}\,]$. By a direct calculation, we obtain
\begin{equation}\label{dyn:ssg}
 \dot{\tau_i}=\Phi_i\tau_i,\quad d_i=\Psi_i \tau_i.
\end{equation}
The system \eqref{dyn:ssg} is called a steady-state generator in \citet{huang2004}, which helps us reject the unwanted disturbance $d_i$. Since the pair $(\Psi_i, \Phi_i)$ is observable, there exists a constant matrix  $G_i$ such that $F_i=\Phi_i+G_i\Psi_i$ is Hurwitz.  To asymptotically reject the disturbance $d_i$, let
\begin{equation}\label{ctr:global-d}
 u_i^d=-\Psi_i\eta_i,\quad \dot{\eta}_i=F_i\eta_i+G_iu_i.
\end{equation}
Inspired by the robust adaptive control law used in \citet{jiang1999robust}, we propose the tracking control with a constant $\lambda>0$ for agent $i$ ($i=1,{\dots},n$):
\begin{align}\label{ctr:global-r}
\begin{split}
u_i^r&=-k_i\rho_i(e_{vi})e_{vi}-\theta_i {\rm sgn}(e_{vi})\\
\dot{k}_i&=-\lambda k_{i}+\rho_i(e_{vi})e_{vi}^2,\quad \dot{\theta}_i=|e_{vi}|
\end{split}
\end{align}
where $e_{vi}=\sum_{j\in \mathcal{V}}a_{ij}(y_i-y_j)$. For simplicity, we take $k_i(0)=\theta_i(0)=0$ (see \citealp{praly2003tac} for a similar setting).

Then the control law for agent $i$ can be written as
\begin{equation}\label{ctr:global}
\begin{split}
u_i&=-\Psi_i\eta_i-k_i\rho_i(e_{vi})e_{vi}-\theta_i {\rm sgn}(e_{vi})\\
\dot{\eta}_i&=F_i\eta_i+G_iu_i\\
\dot{k}_i&=-\lambda k_{i}+\rho_i(e_{vi})e_{vi}^2, \quad \dot{\theta}_i=|e_{vi}|.
\end{split}
\end{equation}

\begin{remark}
Different from most existing internal model (IM) design for multi-agent output regulation \citep[e.g.][]{dong2014asme, su2013scl}, the control (\ref{ctr:global}) contains two parts: the IM design $u_i^d$ for disturbance rejection and the non-smooth design $u_i^r$ to handle the leader's unknown inputs. The gains $k_i$ and $\theta_i$ are designed and updated, independent of $\omega_i(0)$, $v(0)$, and $l$ (refer \citealp{xu2010tcs} and \citealp{tang2014} for similar techniques).
\end{remark}

To select a proper positive function $\rho_i(\cdot)$, we perform a coordinate transformation $\hat \eta_i=\eta_i-\tau_i-G_ie_i$ and the composite system of agent $i$ becomes
\begin{align}\label{dyn:err:im}
\begin{cases}
\dot{\bar z}_i=\bar f_i(\bar z_i, e_i, v, w,\mu)\\
\dot{\hat \eta}_i=F_i\hat\eta_i+F_iG_ie_i- G_i\bar g_i(\bar z_i,e_i, v, w,\mu)\\
\dot{e}_i=\bar g_i(\bar z_i,e_i, v, w,\mu)-\Psi_i\hat \eta_i-\Psi_iG_ie_i+u_i^r\\
\dot{k}_i=-\lambda k_{i}+\rho_i(e_{vi})e_{vi}^2, \quad \dot{\theta}_i=|e_{vi}|.
\end{cases}
\end{align}
Taking $\hat z_i=\text{col}(\bar z_i, \hat \eta_i)$, we have the following result.
\begin{lemma}\label{prop:unknown}
Under Assumption \ref{ass6:unknown}, for any compact subset $\Sigma\subset \mathbb{R}^{\hat n}$, there are smooth Lyapunov functions $V_{\hat z_i}(\cdot)$ satisfying $\hat \alpha_{1i}(||\hat z_i||)\leq V_{\hat z_i}(\hat z_i)\leq \hat \alpha_{2i}(||\hat z_i||)$ for some smooth functions $\hat \alpha_{1i}(\cdot)$ and $\hat \alpha_{2i}(\cdot)\in\mathcal{K}_\infty$ with $i=1,{\dots},n$ such that, for any $(v(t), w(t), \mu)\in \Sigma$ and $\omega_i(t)\in \mathbb{R}^{n_{\omega_i}}$,
\begin{align*}
  \dot{V}_{\hat z_i}(\hat z_i)\mid_{(\ref{dyn:err:im})}\leq -||\hat z_i||^2+\hat \delta_{1i}\hat \gamma_{1i}(e_i)+\hat \delta_{2i}\hat \gamma_{2i}(w)+\hat \delta_{3i}
\end{align*}
where $\hat \gamma_{1i}(\cdot)$ and $\hat \gamma_{2i}(\cdot)$ are known smooth positive definite functions, and $\hat \delta_{1i},\, \hat \delta_{2i},\,\hat \delta_{3i}$ are unknown positive constants.
\end{lemma}
\textbf{Proof.}  From \eqref{err} and $\hat g_i(0, 0, v,\mu)=0$, by Lemma 7.8 in \citet{huang2004}, there are known positive smooth functions $\hat \phi_{1i}(\cdot)$, $\hat \phi_{2i}(\cdot)$ and an unknown constant $\hat c_i>0$ such that
\begin{align}\label{eq:lem1:eq1}
\begin{split}
&||F_iG_ie_i- G_i\bar g_i(\bar z_i,e_i, v, w,\mu)||\\
&\leq \hat c_i (\hat \phi_{1i}(\bar z_i)||\bar z_i||+ \hat \phi_{2i}(e_i)|e_i|+||w||+1).
\end{split}
\end{align}
Based on Assumption \ref{ass6:unknown} and by the changing supply functions technique \citep{sontag1995changing}, for any smooth function $\Delta_i(\bar z_i)>0$, there exists a smooth function $W_{\bar z_i}(\bar z_i)$ satisfying $\bar \alpha_{3i}(||\bar
z_i||)\leq W_{\bar z_i} (\bar z_i)\leq \bar \alpha_{4i}(||\bar z_i||)$
for some smooth functions $\bar \alpha_{3i}(\cdot)$,\,$\bar \alpha_{4i}(\cdot)\in \mathcal{K}_\infty$ such that
  \begin{equation}\label{eq:lem1:eq2}
  \begin{split}
  \dot{W}_{\bar z_i}(\bar z_i)\mid_{(\ref{err})}\leq &-\Delta_i(\bar z_i)||\bar z_i||^2\\
  &+\bar \delta_{3i}\bar \gamma_{3i}(e_i){e_i}^2+\bar \delta_{4i}\bar \gamma_{4i}(w)||w||^2
  \end{split}
  \end{equation}
where $\bar \delta_{3i},\bar \delta_{4i}$ are some unknown positive constants, and $\bar \gamma_{3i}(\cdot), \bar\gamma_{4i}(\cdot)$ are known smooth positive functions.

Let $V_{\hat z_i}(\hat z_i)=\kappa_i W_{\bar z_o}(\bar z_i)+ 2\hat \eta_i^T P_i\hat \eta_i$ with a constant $\kappa_i>0$ to be determined, where $P_i$ is the positive definite solution of the Lyapunov equation $P_i F_i+F_i^TP_i=-I_{n_{p_i}}$.  Clearly, there exist smooth functions $\hat \alpha_{1i}(\cdot)$,\;$\hat \alpha_{2i}(\cdot)\in \mathcal{K}_\infty$ satisfying $\hat \alpha_{1i}(||\hat z_i||)\leq V_{\hat z_i}(\hat z_i)\leq \hat \alpha_{2i}(||\hat z_i||)$.  Using \eqref{eq:lem1:eq1}-\eqref{eq:lem1:eq2} for all $(v(t), w(t),\mu)\in \Sigma$, $\omega_i(t)\in \mathbb{R}^{n_{\omega_i}}$,  we have
\begin{equation*}
\!\begin{split}
  \dot{V}_{\hat z_i}(\hat z_i)\mid_{(\ref{dyn:err:im})} &\leq -[\kappa_i\Delta_i(\bar z_i)-4\hat c_i^2||P_i||^2\hat \phi_{1i}^2(\bar z_i)]||\bar z_i||^2\\
  &-||\hat\eta_i||^2+[\kappa_i\bar \delta_{3i}\bar \gamma_{3i}(e_i)+4\hat c_i^2||P_i||^2\hat \phi_{2i}^2(e_i)]{e_i}^2\\
  &+[\kappa_i\bar \delta_{4i}\bar \gamma_{4i}(w)+4\hat c_i^2||P_i||^2]||w||^2+4\hat c_i^2||P_i||^2
  \end{split}
\end{equation*}
Letting $\kappa_i \geq \max\{1, 4\hat c_i^2||P_i||^2\}$, $\Delta_i(\bar z_i)\geq 1+\hat \phi_{1i}^2(\bar z_i)$, $\hat \delta_{1i}\geq \max \{\kappa_i\bar \delta_{3i}, 4\hat c_i^2||P_i||^2 \}$, $\hat \delta_{2i}\geq \max\{\kappa_i\bar \delta_{4i}, 4\hat c_i^2||P_i||^2\}$, $\hat \delta_{3i}\geq 4\hat c_i^2||P_i||^2$ and  $\hat \gamma_{1i}\geq [\bar \gamma_{3i}(e_i)+\hat \phi_{2i}^2(e_i)]|e_i|^2$, $\hat \gamma_{2i}\geq [\bar \gamma_{4i}(w)+1]||w||^2$ yields the conclusion. \hfill\rule{4pt}{8pt}

It is not hard to check that there is a Filippov solution of  the closed-loop system consisting of \eqref{exosystem}, \eqref{follower}, and \eqref{ctr:global}. Then we present one of our main results.

\begin{theorem}\label{thm3:global:v}
Under Assumptions \ref{ass1:graph}-\ref{ass6:unknown}, there exist smooth positive functions $\rho_i(\cdot)$ for $i=1,{\dots},n$ such that our problem is globally solved by the distributed control (\ref{ctr:global}).
\end{theorem}
\textbf{Proof.}  The proof will be given in two steps.

\emph{Step 1:} We first constructively give smooth positive functions $\rho_i(\cdot)$ for $i=1,{\dots},n$ and show the existence of the trajectory $(z_i, y_i, k_i,\eta_i)$ for $t>0$. Since $\eta_i$-subsystem is linear with a Hurwitz system matrix, we only have to show the boundedness of $(z_i,y_i,k_i)$.  For this purpose, we seek to prove the boundedness of $(\bar z, {e}, \bar k)$, where $\bar k=\text{col}(\bar k_1,{\dots},\bar k_n),\; \bar k_i=k_i-K$ and $K>0$ is a constant to be determined.

From Assumption \ref{ass3:exo}, we can always find a compact subset $\Sigma\subset \mathbb{R}^{\hat n}$ containing $(v(t),w(t),\mu)$ for $t>0$. Based on Assumption \ref{ass6:unknown}, Lemma \ref{prop:unknown} and the changing supply functions technique \citep{sontag1995changing}, for any smooth function $\Delta_i(\hat z_i)>0$, there is a smooth function $W_{\hat z_i}(\hat z_i)$ satisfying $\hat \alpha_{3i}(||\hat z_i||)\leq W_{\hat z_i} (\hat z_i)\leq \hat \alpha_{4i}(||\hat z_i||)$
for some smooth functions $\hat \alpha_{3i}(\cdot)$ and $\hat \alpha_{4i}(\cdot)\in\mathcal{K}_\infty$ such that
\begin{align*}
  \dot{W}_{\hat z_i}(\hat z_i)\mid_{(\ref{dyn:err:im})}&\leq -\Delta_i(\hat z_i)||\hat z_i||^2+\hat \delta_{4i}\hat \gamma_{3i}(e_i){e_i}^2\\
  &+\hat \delta_{5i}\hat \gamma_{4i}(w)||w||^2+\hat \delta_{6i}
\end{align*}
where $\hat \delta_{4i}$, $\hat \delta_{5i}$, $\hat \delta_{6i}$ are some unknown positive constants, and $\hat \gamma_{3i}(\cdot)$, $\hat \gamma_{4i}(\cdot)$ are known smooth positive functions.

Let $W_0(\hat  z)=\sum_{i=1}^n W_{\hat z_i}(\hat z_i), \;V_0(e)={e}^T H{e}.$
By Young's inequality, the set-valued Lie derivative of $V_0$ satisfies
\begin{equation}\label{eq:global:eq1}
\! \dot{V}_0\mid_{(\ref{dyn:err:im})}
\leq 3\sum_{i=1}^ne_{vi}^2 + \sum_{i=1}^n||\Upsilon_{1i}||^2+\Xi_1+ \bigcap_{\xi\in \partial V_0}\xi^T u^r
\end{equation}
where $\Upsilon_{1i}\triangleq \hat g_i(\bar z_i,e_i, v,\mu)-\Psi_i\hat \eta_i-\Psi_iG_ie_i$ and $\Xi_1\triangleq \sum_{i=1}^n(\bar b_2^2||w||^2+||{\bf u}_i(v,\mu)||^2)$.
Since $\hat g_i(0,0,v,\mu)=0$ and $v$ is bounded, by Lemma 7.8 in \citet{huang2004}, there exist known positive smooth functions $\phi_{1i}(\cdot), \phi_{2i}(\cdot)$ with an unknown positive constant $c_i$ such that
\begin{equation}\label{eq:global:eq2}
  ||\Upsilon_{1i}||^2\leq c_i(\phi_{1i}(\hat z_i)||\hat z_i||^2+\phi_{2i}(e_i)||e_i||^2).
\end{equation}
Note that $$\sum_{i=1}^n[\hat \delta_{4i}\hat \gamma_{3i}(e_i)+c_i \phi_{2i}(e_i)]{e_i}^2\leq \sum_{i=1}^n  \delta_{i}\hat\phi_{i}(e_i){e_i}^2$$
with $\delta_i=\hat \delta_{4i}+c_i$ and a positive smooth function $\hat \phi_i(e_i)\geq
\max\{\hat \gamma_{3i}(e_i),\phi_{2i}(e_i)\}$.  Recalling $e_v=He$, by
similar arguments, we obtain
\begin{equation}\label{eq:global:eq4}
\sum_{i=1}^n[\hat \delta_{4i}\hat\gamma_{3i}(e_i)+ c_i\phi_{2i}(e_i)]{e_i}^2\leq \sum_{i=1}^n\delta_i \tilde\phi_i(e_{vi})e_{vi}^2
\end{equation}
for some known positive smooth functions $\tilde\phi_i(\cdot)$.

Construct a positive definite and radially unbounded Lyapunov function candidate as
$
V(\hat z, {e}, \bar k)=W_0(\hat z)+V_0(e)+\sum_{i=1}^n\bar k_i^2.
$
Clearly, $\theta_i\geq0$, and by (\ref{eq:global:eq1})-(\ref{eq:global:eq4}), it follows
\begin{align*}
\dot{V}\mid_{(\ref{dyn:err:im})}\leq & \sum_{i=1}^n-[\Delta_i(\bar z_i)- c_i\phi_{1i}(\hat z_i)]||\hat z_i||^2-\sum_{i=1}^n\lambda \bar k_i^2 \\
  &- \sum_{i=1}^n[2K\rho(e_{vi})-3-\delta_i \tilde\phi_i(e_{vi})]e_{vi}^2+\Xi_2
\end{align*}
where \[\Xi_2\triangleq \Xi_1+\sum_{i=1}^n \hat \delta_{5i}\hat \gamma_{4i}(w)||w||^2+\sum_{i=1}^n\hat \delta_{6i}+\sum_{i=1}^n\lambda K^2.\]
For $i=1,{\dots},n$, we take $\Delta_i, \rho_i$ and $K$ with $\Delta_i(\hat z_i)=-c_i\phi_{1i}(\hat z_i)+1,\, \rho_i(e_{vi})=\tilde\phi_i(e_{vi})+1$, $K>\frac{1}{2}(\delta_i+4)$. Then
\begin{equation*}
  \dot{V}\mid_{(\ref{dyn:err:im})}\leq-\sum_{i=1}^n||\hat z_i||^2 - \sum_{i=1}^n e_{vi}^2 -\sum_{i=1}^n\lambda \bar k_i^2+\Xi_2.
\end{equation*}
Similar to the analysis of Theorem 4.14 in \citet{khalil2002}, for any $\lambda>0$, the uniform boundedness of $\hat z$, $e$ and $\bar k$ with any solution of the closed-loop system is obtained, which implies the boundedness of $z_i, y_i$, and $k_i$.

\emph{Step 2:} Let us check ${e}$ under the controller (\ref{ctr:global}).  Take a positive definite function $V_1({e},\bar\theta)=V_0(e)+\sum_{i=1}^n \bar \theta_i^2$, where $\bar \theta_i=\theta_i-\Theta$ with $\Theta>0$ to be specified later. Then
\begin{align*}
  \dot{V}_1\mid_{(\ref{dyn:err:im})} &=2\sum_{i=1}^ne_{vi}[\Upsilon_{2i}-\theta_i {\rm sgn}(e_{vi})]+2\sum_{i=1}^n {\bar \theta_i}\dot{\theta}_i\\
  &\leq 2\sum_{i=1}^ne_{vi}[\Upsilon_{2i}-\Theta {\rm sgn}(e_{vi})]
\end{align*}
where $\Upsilon_{2i}\triangleq \bar g_i(\bar z_i, e_i, v, w,\mu)-\Psi_i\hat \eta_i-\Psi_iG_ie_i-k_i\rho_i(e_{vi})e_{vi}$.

From the boundedness of $\hat z, e, v, \omega,\mu$ and the continuity of related functions, there is a positive constant $\hat M$ such that $||\Upsilon_{2i}||\leq \hat M$ for all $i$.
Letting $\Theta > \hat M+1$ gives
\begin{align} \label{thm3:global:v:integ}
  \dot{V}_1\mid_{(\ref{dyn:err:im})} &\leq -2\sum_{i=1}^n(\Theta -\hat M)|e_{vi}|\leq -2\sum_{i=1}^n|e_{vi}|.
\end{align}
Because $V_1$ is radially unbounded, $\theta_i$ is bounded. From the boundedness of $e$ showed above, $\dot{e}_i$ is also bounded by (\ref{ctr:global}) and (\ref{dyn:err:im}), and then we can obtain the uniform continuity of $V_2 \triangleq 2\sum_{i=1}^n|e_{vi}|$ with respect to
time $t$. Integrating (\ref{thm3:global:v:integ}) from $0$ to $t$
and taking $t\to \infty$, we have
$$
\int_{0}^{\infty} V_2({e}(t)){\mathrm d}t\leq V_1(0).
$$
Recalling the Barbalat's lemma \citep{khalil2002}, $V_2({e}(t))\to 0$
when $t\to \infty$, and hence $e_i$
converges to zero as $t\to\infty$ for $i=1,\dots,n$. Thus,
the conclusion is obtained. \hfill\rule{4pt}{8pt}

\subsection{Special Case}

In some cases, we may know the domain of the initial condition $(z_i(0),y_i(0))$, $v(0)$, $\omega_i(0)$ and the bound of the unknown input and the uncertain parameter $\mu$.  Of course, we can still use the proposed control law (\ref{ctr:global}), but here we construct a simpler control law to solve the problem based on the additional information.  To this end, it is reasonable to introduce a new assumption to replace Assumption \ref{ass6:unknown}.

\begin{assumption}\label{ass5:zero-known}
Given any compact subset $\Sigma\subset \mathbb{R}^{\hat n}$, there exist smooth Lyapunov functions $V_{\bar z_i}(\cdot)$ satisfying $\alpha_{1i}(||\bar z_i||)\leq V_{\bar z_i}(\bar z_i)\leq \alpha_{2i}(||\bar z_i||)$ for some smooth functions $\alpha_{1i}(\cdot)$ and $\alpha_{2i}(\cdot)\in \mathcal{K}_\infty$ (for $i=1,\dots,n$) such that, for any $(v(t), w(t), \mu)\in \Sigma$,
  \begin{equation*}
  \dot{V}_{\bar z_i}\mid_{(\ref{err})}\leq -\alpha_i(||\bar z_i||)+\gamma_{1i}(e_i)+\gamma_{2i}(w)
  \end{equation*}
where $\gamma_{1i}(\cdot)$ and $\gamma_{2i}(\cdot)$ are known smooth positive definite functions, and $\alpha_i(\cdot)$ is a known class $\mathcal{K}_\infty$ function satisfying $\limsup_{s \to 0+}(\alpha_i^{-1}(s^2)/s)\leq \infty$.
\end{assumption}

Then a new simple controller is proposed in this case:
\begin{align}\label{ctr:semi}
\begin{split}
u_i&=-\Psi_i\eta_i-\rho_i(e_{vi})e_{vi}-\gamma_i {\rm sgn}(e_{vi})\\
\dot{\eta}_i&=F_i\eta_i+G_iu_i.
\end{split}
\end{align}
Performing a transformation $\hat \eta_i=\eta_i-\tau_i-G_i e_i$ gives
\begin{align}\label{dyn:err:im-static}
\begin{cases}
\dot{\bar z}_i=\bar f_i(\bar z_i, e_i, v, w,\mu)\\
\dot{\hat \eta}_i=F_i\hat\eta_i+F_iG_ie_i-G_i\bar g_i(\bar z_i,e_i, v, w,\mu)\\
\dot{e}_i=\bar g_i(\bar z_i,e_i, v, w,\mu)-\Psi_i\hat \eta_i-\Psi_iG_ie_i+u_i^r.
\end{cases}
\end{align}
Denote $\hat z_i=\text{col}(\bar z_i, \hat \eta_i)$, and the next lemma can be proved in a similar way given in the last subsection.
\begin{lemma}\label{prop:known}
Under Assumption \ref{ass5:zero-known}, for any given positive constant $M$ and for $i=1,{\dots},n$, there is a positive constant $\bar M$ and a smooth Lyapunov function $V_{\hat z_i}(\cdot)$ satisfying $\hat \alpha_{1i}(||\hat z_i||)\leq V_{\hat z_i}(\hat z_i)\leq \hat \alpha_{2i}(||\hat z_i||)$ for some smooth functions $\hat \alpha_{1i}(\cdot)$, $\hat \alpha_{2i}(\cdot)\in\mathcal{K}_\infty$, such that, for all $(v(0),\omega_i(0),\mu)\in \mathbb{R}_M^{\hat n}$, $w(t)\in \mathbb{R}_M^{n_w}$ and $\eta_i(0)\in \mathbb{R}_{\bar M}^{n_{\eta_i}}$,
  \begin{align*}
  \dot{V}_{\hat z_i}(\hat z_i)\mid_{(\ref{dyn:err:im-static})}\leq -||\hat z_i||^2+ \hat \gamma_{1i}(e_i)+\hat \gamma_{2i}(w)+\hat \delta
  \end{align*}
where $\hat \gamma_{1i}(\cdot), \hat \gamma_{2i}(\cdot)$ are known smooth positive definite functions and $\hat \delta$ is a known positive constant.
\end{lemma}

Then we show how the control (\ref{ctr:semi}) solves our problem.

\begin{theorem}\label{thm1:semi}
Under Assumptions \ref{ass1:graph}-\ref{ass4:re} and \ref{ass5:zero-known}, our problem can be semi-globally solved by the distributed control \eqref{ctr:semi}.
\end{theorem}
\textbf{Proof.} The proof is similar to that of Theorem \ref{thm3:global:v}.

\emph{Step 1:} We first prove the boundedness of $(\hat z, {e})$.
Based on Assumption \ref{ass5:zero-known} and Lemma \ref{prop:known}, we apply the changing supply functions technique and obtain that, for any given $M>0$ and any smooth function $\Delta_i(\hat z_i)$, there is a positive constant $\bar M$ and a smooth function $W_{\hat z_i}(\hat  z_i)$ satisfying $\hat\alpha_{3i}(||\hat z_i||)\leq W_{\hat z_i}(\hat z_i)\leq
\hat\alpha_{4i}(||\hat z_i||)$ for some smooth functions $\hat\alpha_{3i}(\cdot),\;\hat\alpha_{4i}(\cdot)\in\mathcal{K}_\infty$, such that, for any $(v(0),\omega_i(0),\mu)\in \mathbb{R}_M^{\hat n}$, $w(t)\in \mathbb{R}_M^{n_w}$ and $\eta_i(0)\in \mathbb{R}_{\bar M}^{n_{\eta_i}}$,
  \begin{align}\label{eq:semi:eq1}
  \begin{split}
  \dot{W}_{\hat z_i}\mid_{(\ref{dyn:err:im-static})}\leq&-\Delta_i(\hat z_i)||\hat
  z_i||^2\\
  &+\hat\gamma_{3i}(e_i){e_i}^2+\hat\gamma_{4i}(w)||w||^2+\hat \delta_1
  \end{split}
  \end{align}
where $\hat\gamma_{3i}(\cdot),\hat\gamma_{4i}(\cdot)$ are known smooth positive functions and $\hat \delta_1$ is a known positive constant.

Consider a positive definite and radially unbounded function $\bar V(\hat z, {e})=W_0(\hat z)+V_0(e)$. For a given $M>0$,  by similar arguments in Theorem \ref{thm3:global:v}, there exist known smooth positive functions $\phi_{1i},\phi_{2i}, \tilde\phi_i$ such that
\begin{align}\label{eq:semi:eq2}
\begin{split}
  ||\Upsilon_{3i}||^2 \leq\phi_{1i}(\hat z_i)||\hat z_i||^2+\phi_{2i}(e_i)||e_i||^2\\
  \sum_{i=1}^n[\hat \gamma_{3i}(e_i)+ \phi_{2i}(e_i)]{e_i}^2\leq\sum_{i=1}^n \tilde\phi_i(e_{vi})e_{vi}^2
  \end{split}
\end{align}
where $\Upsilon_{3i}\triangleq\hat g_i(\bar z_i,e_i, v, \mu)-\Psi_i\hat \eta_i-\Psi_iG_ie_i$. By (\ref{eq:global:eq1}), (\ref{eq:semi:eq1}), (\ref{eq:semi:eq2}) and $\gamma_i>0$,  letting $\Delta_i(\hat z_i)=-\phi_{1i}(\hat z_i)+1, \rho_i(e_{vi})=\frac{1}{2}(\tilde\phi_i(e_{vi})+4)$ yields $ \dot{\bar V}\mid_{(\ref{dyn:err:im-static})} \leq -||\hat z||^2-||e_v||^2+\Xi_3$,
where $$\Xi_3\triangleq\sum_{i=1}^n[\hat \gamma_{4i}(w)+\bar  b_2^2]||w||^2+n\hat\delta_1+\sum_{i=1}^n||{\bf u}_i(v,\mu)||^2.$$

Again by similar techniques used in Theorem 4.14 in \citet{khalil2002}, we obtain the uniform
boundedness of the trajectory $(\hat z,e)$ and a positive constant
$\hat M$, only depending on $M$ and $\bar M$, satisfying $
||\Upsilon_{4i}||\leq \hat M$ with $\Upsilon_{4i}\triangleq\bar g_i(\bar z_i,e_i, v, w,\mu)-\Psi_i\hat \eta_i-\Psi_iG_ie_i-\rho_i(e_{vi})e_{vi}$.

\emph{Step 2:} Check ${e}$ by considering the set-valued Lie derivative of $V_0({e})={e}^T H{e}$. Taking $\gamma_i \geq \hat M+1$ gives
\begin{align*}
  \dot{ V}_0\mid_{(\ref{dyn:err:im-static})}&=2\sum_{i=1}^ne_{vi}^T[\Upsilon_{4i}-\gamma_i {\rm sgn}(e_{vi})]\\
  &\leq -2\sum_{i=1}^n(\gamma_i -\hat M)|e_{vi}|\leq -2\sum_{i=1}^n|e_{vi}|.
\end{align*}
Thus, $\lim_{t\to +\infty}e(t)=0$, which implies the conclusion. \hfill\rule{4pt}{8pt}

\section{Simulations}

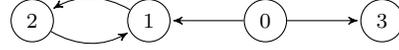
\begin{figure}
\centering
\begin{tikzpicture}[shorten >=1pt, node distance=1.55 cm, >=stealth',
every state/.style ={circle, minimum width=0.25cm, minimum height=0.25cm}]
\node[align=center,state](node2) {\scriptsize 2};
\node[align=center,state](node1)[right of=node2]{\scriptsize 1};
\node[align=center,state](node0)[right of=node1]{\scriptsize 0};
\node[align=center,state](node3)[right of=node0]{\scriptsize 3};
\path[->]   (node0) edge (node1)
            (node1) edge [bend right] (node2)
            (node2) edge [bend right] (node1)
            (node0) edge (node3)
            ;
\end{tikzpicture}
\caption{Interaction graph $\mathcal G$ in our
example.}\label{fig:graph}
\end{figure}

\begin{figure}
  \centering
  \subfigure[]{
    \includegraphics[width=0.22\textwidth]{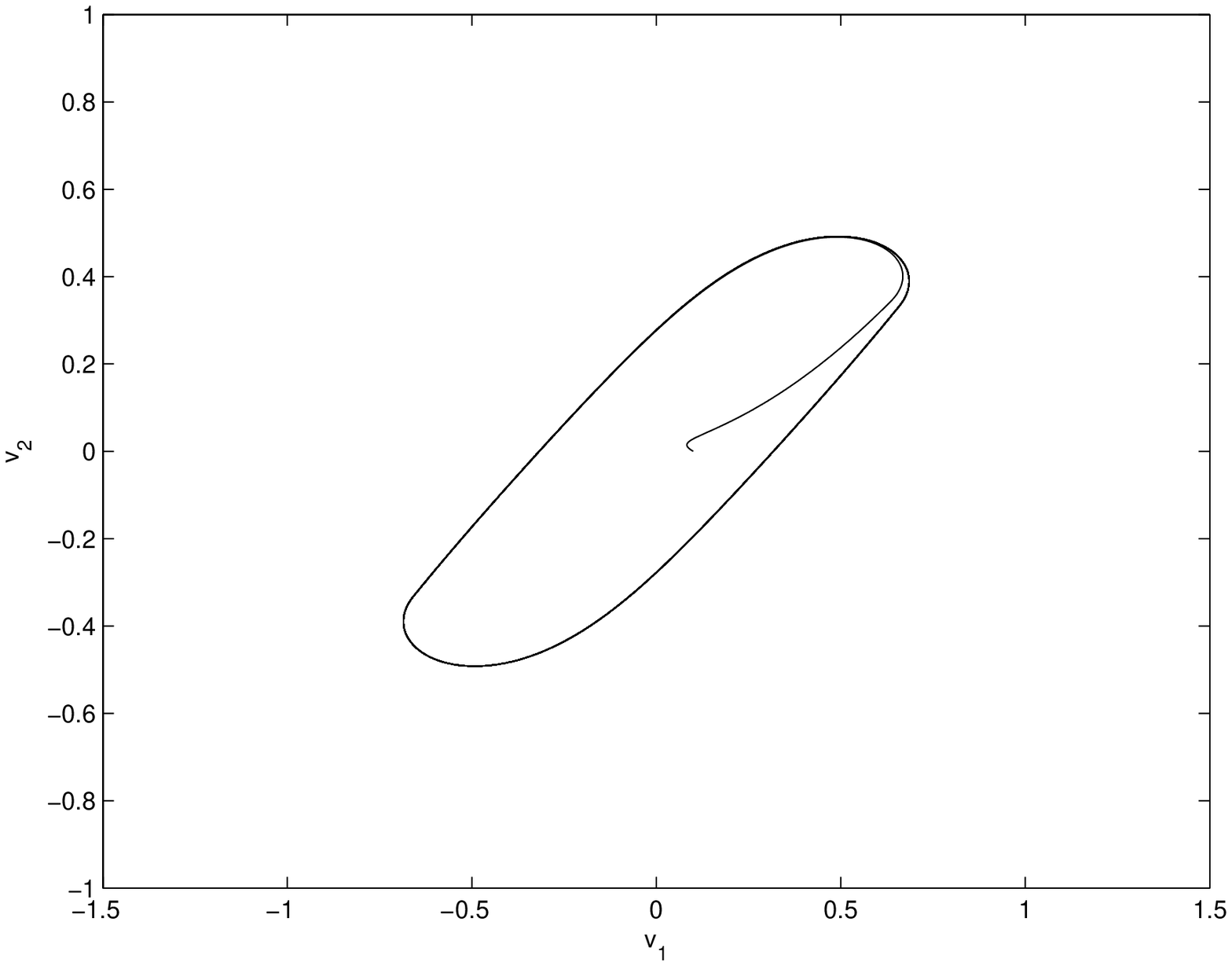}
  }
  \subfigure[]{
    \includegraphics[width=0.22\textwidth]{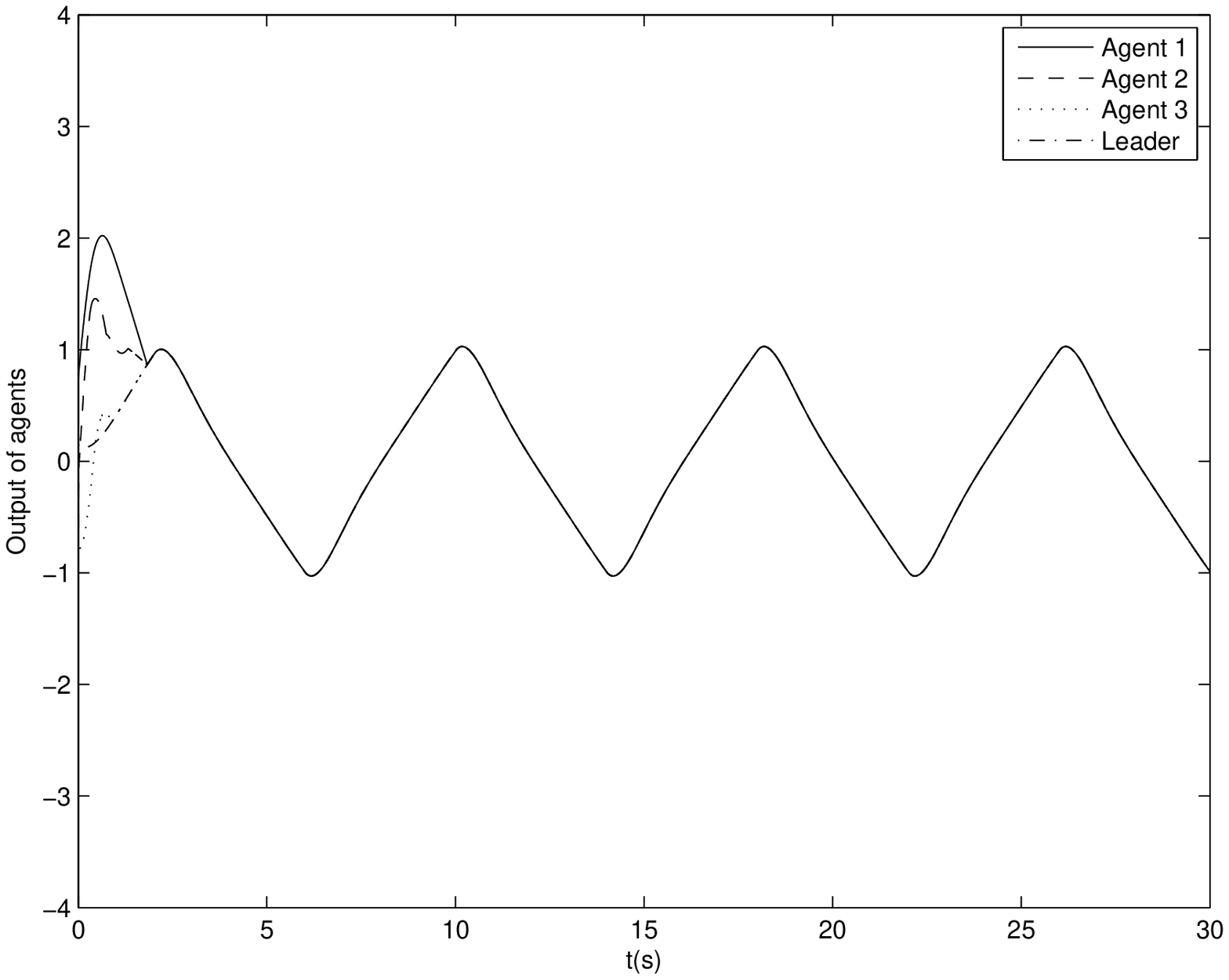}
    }
\caption{Trajectory of the leader and performance of the control law.}\label{fig:simu}
\end{figure}

To illustrate our control design, we consider a multi-agent system with three non-identical followers in the form of high-order FitzHugh-Nagumo dynamics \citep{xu2010tcs} as follows:
\begin{equation*}
    \begin{cases}
        \dot{x}_{1i}=x_{1i}-\frac{1}{3}x_{1i}^3-x_{2i}+x_{3i}+d_i(t)+b_iu_i\\
        \dot{x}_{2i}=x_{1i}+c_{1i}-x_{2i}\\
        \dot{x}_{3i}=-x_{1i}+c_{2i}-x_{3i}, \qquad i=1,2,3
    \end{cases}
\end{equation*}
where $c_{1i},\,c_{2i},$ and $b_i$ are positive constants. The local disturbances are generated by (\ref{disturbance}) with $D_1=1+\mu_1, \; S_1=0, \; D_2=[1+\mu_2,0],\; S_2=[0,1; 0,0],\; D_3=[1+\mu_3,0],\; S_3=[0, 1; -1,0]$. We aim to make $x_{1i}$ follow a reference $y_0=(1+\mu_v)v_1$ generated by a leader \citep{grasman1987} with the input $p(t)$ unknown to the followers as follows:
\begin{align*}
  \begin{cases}
    \dot{v}_{1}=\epsilon_0(-2v_1-v_2^3-v_2)+p(t)\\
    \dot{v}_{2}=v_1-v_{2}
    \end{cases}
\end{align*}
where $p(t)=\frac{2A}{T}\left(t-T\lfloor \frac{t}{T}+\frac{1}{2}\rfloor\right)(-1)^{\lfloor \frac{t}{T}+\frac{1}{2}\rfloor}$ is a triangle wave signal with period $2T$ and amplitude $A$, and $\lfloor x \rfloor$ is the largest integer not greater than $x$. Denote $\mu \triangleq \mbox{col}(\mu_v,\mu_x,\mu_1,\mu_2,\mu_3)$ as the uncertain parameter vector.  Figure \ref{fig:graph} describes this multi-agent interaction topology with $a_{ij}=1$, and Figure \ref{fig:simu}(a) depicts the reference trajectory with $v_1(0)=0.1$ and $v_2(0)=0$.

Without knowing the exact form of $p(t)$, the formulation in  \citet{su2013scl} even with a nonlinear exosystem \citep{chen2005auto} fails to solve this problem for our multi-agent system. Nevertheless, the problem is solvable based on the formulation in Section 2. In fact, Assumption \ref{ass3:exo} is verified by a Lyapunov function $ V = v_1^2 + \frac{\epsilon_0}{2}v_2^4+\epsilon_0 v_2^2$. Let $z_{1i}=x_{2i}-c_{1i}, z_{2i}=x_{3i}-c_{2i},\, y_i=x_{1i}$, and the plant is of the form \eqref{follower} satisfying Assumption \ref{ass4:re} with ${\bf z}_{1i}=(1+\mu_v)v_2,\,{\bf z}_{2i}=-(1+\mu_v)v_2$. Assumption \ref{ass6:unknown} also holds with $V_{\bar z_i}(s)=\frac{1}{2}s^2, \; \alpha_i(s)=\frac{1}{2}s^2,\; \gamma_{1i}(s)=\gamma_{2i}(s)=s^2.$  By Theorem \ref{thm3:global:v}, the control with $\rho_i(s)=s^6+1$ can solve this problem. To reduce the unfavorable chattering in simulations, the sign function in the proposed control \eqref{ctr:global} can be replaced by a saturation function as follows:
 \[
{\rm sat}_\epsilon(x)=
\begin{cases}
x/\epsilon, &\text{ if } |x|\leq\epsilon;\\
{\rm sgn}(x/\epsilon),& \text{ if } |x|>\epsilon.
\end{cases}
\]
With selected matrix pairs $(F_i,\,G_i)$, Figure \ref{fig:simu}(b) shows the performance of the controller with $A=2,\, T=4$,\, $\epsilon=10^{-3}$,\,$c_{1i}=c_{2i}=b_i=i$ and $\mu\in [-1,1]^5$.

\section{Conclusions}
In this paper, a distributed output regulation problem was formulated for a class of uncertain heterogeneous nonlinear multi-agent systems to deal with local disturbances and an unknown-input leader. Based on changing supply functions and adaptive techniques, distributed control laws incorporating local internal models were constructed. The semi-global and global results were obtained in two different cases.

\end{document}